\documentstyle[12pt,leqno]{article}

\newtheorem{thm}{Theorem}[section]
\newtheorem{Theorem}{Theorem}[section]

\newtheorem{lem}[thm]{Lemma}
\newtheorem{pro}[thm]{Proposition}
\newcommand{\sthm}{\begin{Theorem}}         
\newcommand{\ethm}{\end{Theorem}}           
\newtheorem{Corollary}[Theorem]{Corollary}
\newcommand{\scor}{\begin{Corollary}}       
\newcommand{\ecor}{\end{Corollary}}         
\newcommand{\pf}{ \par \vspace{1ex} \noindent {\sc Proof.} \hspace{2mm}}

\begin{document}
\title{Lagrangian $H$-umbilical submanifolds in quaternion Euclidean spaces}
\author{Yun Myung Oh\\
ohy@math.msu.edu}
\date{ }
\maketitle
\begin{abstract}
In \cite{c97}, B.Y. Chen proved that the Lagrangian $H$-umbilical
submanifolds in complex Euclidean space $\bf C^{n}$ are Lagrangian
pseudo-spheres and complex extensors of the unit hypersphere of
$\bf E^{n}$, except the flat ones. Similar to this, we can define
the Lagrangian $H$-umbilical submanifold in quaternion space
forms. The main purpose of this paper is to classify the
Lagrangian $H$-umbilical submanifolds in quaternion Euclidean
space $\bf H^{n}$.
\end{abstract}
\section{Introduction}
It has been known that there is no totally umbilical Lagrangian
submanifolds in complex-space-forms except the totally geodesic
ones. It was natural to look for the simplest Lagrangian
submanifold next to the totally geodesic ones in
complex-space-form. To do so, B.Y. Chen introduced the notion of
the Lagrangian $H$-umbilical submanifold \cite{c97b} and also, in
\cite{c97}, he obtained the classification theorems for Lagrangian
$H$-umbilical submanifolds in complex
Euclidean space $\bf C^{n}$.\\
Similar to the above case, it also has been known that there exist
no totally umbilical Lagrangian submanifold in
quaternion-space-form $\widetilde{M}^{n}(4c)$ except the totally
geodesic ones. In order to find the next simplest case in
quaternion-space-form $\widetilde{M}^{n}(4c)$, we also introduce
the notion of Lagrangian $H$-umbilical submanifold. By a
Lagrangian $H$-umbilical submanifold $M^{n}$ in a quaternion
manifold $\widetilde{M}^{n}$ we mean a non-totally geodesic
Lagrangian submanifold whose second fundamental form is given by
\begin{equation}
\begin{array}{lll}
h(e_{1},e_{1}) & =&  \lambda_{1} Ie_{1} + \lambda_{2} Je_{1} + \lambda_{3} Ke_{1},\\
h(e_{2},e_{2}) & =&  .... = h(e_{n},e_{n}) = \mu_{1} Ie_{1} + \mu_{2} Je_{1} + \mu_{3} Ke_{1},\\
h(e_{1},e_{j}) & =&  \mu_{1} Ie_{j} + \mu_{2} Je_{j} + \mu_{3} Ke_{j},\hskip1pc j=2,...,n\\
h(e_{j},e_{k}) & =& 0,\hskip1pc j \neq k, j,k=2,...,n
\end{array}\label{eq:1}
\end{equation}
for some functions $\lambda_{1},\lambda_{2},\lambda_{3},\mu_{1},\mu_{2}$ and
$\mu_{3}$ with respect to an orthonormal local frame fields.\\
According to the above condition, the mean curvature vector
$H$ is given by $H=H_{1}+H_{2}+H_{3}$, where
$H_{1}=\gamma_{1}Ie_{1}$,$H_{2}=\gamma_{2}Je_{1}$, $H_{3}=\gamma_{3}Ke_{1}$,
and $\displaystyle{\gamma_{i}=\frac{\lambda_{i}+(n-1)\mu_{i}}{n}(i=1,2,3)}$.
The condition (\ref{eq:1}) is equivalent to
\begin{eqnarray*}
h(X,Y) & = & \alpha_{1}<IX,H><IY,H>H_{1}+\alpha_{2}<JX,H><JY,H>H_{2}\\
       &   & +\alpha_{3}<KX,H><KY,H>H_{3})\\
       &   & +\beta_{1}(<X,Y>H_{1}+<IY,H>IX+<IX,H>IY)\\
       &   & +\beta_{2}(<X,Y>H_{2}+<JY,H>JX+<JX,H>JY)\\
       &   & +\beta_{3}(<X,Y>H_{3}+<KY,H>KX+<KX,H>KY)
\end{eqnarray*}
for any tangent vectors $X,Y$, where
$\alpha_{i}=\displaystyle{\frac{\lambda_{i}-3\mu_{i}}{\gamma_{i}^{3}}}$
and $\displaystyle{\beta_{i}=\frac{\mu_{i}}{\gamma_{i}}}$ for $i=1,2,3$.
Clearly, a non-minimal Lagrangian $H$-umbilical submanifold has the
shape operator $A_{H}$ at $H$
with two eigenvalues $\lambda$ and $\mu$,
where $\lambda=\sum_{i=1}^{n}\lambda_{i}\gamma_{i}$ and
$\mu=\sum_{i=1}^{n}\mu_{i}\gamma_{i}$ with respect to some orthonormal frame fields.\\
On the other hand, it also satisfies
$<h(X,Y),\phi_{i}Z>=<h(Z,Y),\phi_{i}X>$,
where $\phi_{i}$ is one of the element in $\{I,J,K\}$ and $X,Y,Z$ are tangent vectors  to $M^{n}$.
Using this property, we can say that Lagrangian $H$-umbilical
submanifolds are the simplest Lagrangian submanifolds next to totally geodesic
submanifolds in quaternion Euclidean space.\\
The main purpose of this paper is to classify the Lagrangian $H$-umbilical
submanifolds in quaternion Euclidean space.
\section{Preliminaries \cite{is}}
Let $\widetilde{M}^{n}$ be a $4n$-dimensional Riemannian manifold
with metric $g$. $\widetilde{M}^{n}$ is called a quaternion
manifold if there exists a 3-dimensional vector space $V$ of
tensors of type (1,1) with local basis of almost Hermitian
structure $I,J$ and $K$ such that
\\[0.1in]
$\left.\begin{tabular}{l}
(a)$IJ=-JI=K, JK=-KJ=I, KI=-IK=J, I^{2}=J^{2}=K^{2}=-1,$\\
(b) for any local cross-section $\varphi$ of $V$,$\widetilde{\nabla}_{X}\varphi$
is also a cross section of $V$,\\
\end{tabular}\right.$\\
when $X$ is an arbitrary vector field on $\widetilde{M}^{n}$ and $\widetilde{\nabla}$
the Riemannian connection on $\widetilde{M}^{n}$.
\\[0.1in]
Condition (b) is equivalent to the following condition:\\
$(b')$ there exist local 1-forms $p,q$ and $r$ such that\\
$$\left.\begin{array}{ccccc}
\widetilde{\nabla}_{X}I & = & & r(X)J-q(X)K,\\
\widetilde{\nabla}_{X}J & = & -r(X)I & \mbox{} + p(X)K,\\
\widetilde{\nabla}_{X}K & = & q(X)I & - p(X)J
\end{array}\right.$$
Let $X$ be a unit vector on $\widetilde{M}^{n}$. Then $X,IX,JX$, and $KX$ form
an orthonormal frame on $\widetilde{M}^{n}$. We denote by $Q(X)$ the 4-plane
spanned by them. For any two orthonormal vectors $X,Y$ on $\widetilde{M}^{n}$,
if $Q(X)$ and $Q(Y)$ are orthogonal, the plane $\pi(X,Y)$ spanned by $X,Y$ is
called a totally real plane. Any 2-plane in a $Q(X)$ is called a quaternion plane.
The sectional curvature of a quaternion plane $\pi$ is called the quaternion
sectional curvature of $\pi$. A quaternion manifold is a quaternion-space-form
if its quaternion sectional curvatures are equal to a constant $4c$. We denote
such a $4n$-dimensional quaternion-space-form by $\widetilde{M}^{n}(4c)$.
\\[0.2in]
It is well known that a quaternion manifold $\widetilde{M}^{n}$ is a
quaternion-space-form if and only if its curvature tensor $\widetilde{R}$ is
of the following form:
\begin{eqnarray}\label{eq:2}
\widetilde{R}(X,Y)Z & = & c\{g(Y,Z)X - g(X,Z)Y\\\nonumber
                &   & \mbox{}+ g(IY,Z)IX - g(IX,Z)IY + 2g(X,IY)IZ\\\nonumber
                &   & \mbox{}+ g(JY,Z)JX - g(JX,Z)JY + 2g(X,JY)JZ\\\nonumber
                &   & \mbox{}+ g(KY,Z)KX - g(KX,Z)KY + 2g(X,KY)KZ \}\nonumber
\end{eqnarray}
for tangent vectors $X,Y$ and $Z$ on $\widetilde{M}^{n}$.
\\[0.1in]
Let $M$ be an $n$-dimensional Riemannian manifold and $x:M\rightarrow \widetilde{M}^{n}(4c)$
be an isometric immersion of $M$ into a quaternion-space-form $\widetilde{M}^{n}(4c)$.
We call $M$ a Lagrangian submanifold or a totally real submanifolds of $\widetilde{M}^{n}(4c)$
if each 2-plane of $M$ is mapped by $x$ into a totally real plane in $\widetilde{M}^{n}(4c)$.
Consequently if $M$ is a Lagrangian submanifold of $\widetilde{M}^{n}(4c)$ then
$\phi(TM)\subset T^{\bot}M$ for $\phi=I,J,$ or $K$, $T^{\bot}M$ is the normal bundle of
$M$ in $\widetilde{M}^{n}(4c)$.
\\[0.1in]
For any orthonormal vectors $X,Y$ in $Y$, $\pi(X,Y)$ is totally real in
$\widetilde{M}^{n}(4c)$, $Q(X)$ and $Q(Y)$ are orthogonal and $g(X,\varphi Y)=g(\phi X,Y)=0$
for $\varphi,\phi=I,J$ or $K$. By (\ref{eq:2}) we have
$$\widetilde{R}(X,Y)Z  =  c\{g(Y,Z)X - g(X,Z)Y\},\mbox{   for  }X,Y,Z\in TN$$
If we denote the Levi-Civita connections of $M$ and $\widetilde{M}^{n}(4c)$ by
$\nabla$ and $\tilde{\nabla}$, respectively, the formulas of Gauss and Weingarten
are respectively given by
$$\left.\begin{array}{l}
\tilde{\nabla}_{X}Y = \nabla_{X}Y + h(X,Y),\\
\tilde{\nabla}_{X}\zeta = -A_{\zeta}X + D_{X}\zeta,
\end{array}\right.$$
for tangent vector fields $X$, $Y$ and normal vector field $\zeta$,
where $D$ is the normal connection. The second fundamental form $h$
is related to the shape operator
$A_{\zeta}$ by $\langle h(X,Y),\zeta \rangle = \langle A_{\zeta}X,Y \rangle$.
The mean curvature vector $H$ of $M$ is defined by $H = \frac{1}{n}$ trace $h$.
If we denote the curvature tensors of $\nabla$ and $D$ by $R$ and $R^{D}$,
then the equations of Gauss, Codazzi and Ricci are given by
$$\left.\begin{array}{cll}
\langle R(X,Y)Z,W \rangle & = & \langle h(Y,Z),h(X,W) \rangle - \langle h(X,Z),h(Y,W) \rangle\\
 & & + c(\langle X,W \rangle\langle Y,Z \rangle - \langle X,Z \rangle\langle Y,W \rangle),\\
(\nabla h)(X,Y,Z) & = & (\nabla h)(Y,X,Z),\\
\langle R^{D}(X,Y)Z,W \rangle & = & \langle [A_{JZ},A_{JW}]X,Y \rangle +
c(\langle X,W \rangle\langle Y,Z \rangle - \langle X,Z \rangle\langle Y,W \rangle),
\end{array}\right.$$
where $X,Y,Z,W$ are tangent vector fields and $\eta$ and $\zeta$ are normal vector fields to $M$.
\vskip.2pc
Finally, we recall a definition of warped product \cite{c81}. Let
$N_{1},N_{2}$ be two Riemannian manifolds with Riemannian metrics
$g_{1},g_{2}$, respectively and $f$ a positive function on
$N_{1}$. Then the metric $g=g_{1}+f^{2}g_{2}$ is called a warped
product metric on $N_{1}\times N_{2}$. The manifold $N_{1}\times
N_{2}$ with the warped product metric $g=g_{1}+f^{2}g_{2}$ is
called a warped product manifold. The function $f$ is called the
warping function of the warped product manifold.
\section{Quaternion extensors}
In this section we are going to investigate the geometry of quaternion extensors.
First of all, we have the following definition.\\
Let $G:M^{n-1}\rightarrow \bf E ^{m}$ be an isometric immersion of a Riemannian
(n-1)-manifold into Euclidean m-space $\bf E^{m}$ and $F:I\rightarrow \bf H$ a
unit speed curve in the quaternion plane. Consider the following extension $\phi$ given by
$$\phi=F\otimes G:I\times M^{n-1}\rightarrow \bf H\otimes \bf E^{m}=\bf H^{m},$$
where $\phi=F\otimes G$ is the tensor product immersion of $F$ and $G$ defined by
$$(F\otimes G)(s,p)=F(s)\otimes G(p);\mbox{   }s\in I,p\in M^{n-1}.$$
We call such an extension $\phi=F\otimes G$ a quaternion extensor of G via F.\\
An immersion $f:N\rightarrow \bf E^{m}$ is called spherical
(respectively, unit spherical) if $N$ is immersed into a
hypersphere (respectively, unit hypersphere) of $\bf E^{m}$
centered at the origin. The quaternion extensor $\phi:F\otimes
G:I\times M^{n-1}\rightarrow \bf H^{m}$ is called $F$-isometric
if, for each $p\in M^{n-1}$, the immersion $F\otimes
G(p):I\rightarrow {\bf {H^{m}}}:s\mapsto F(s)\otimes G(p)$ is
isometric. Similarly, the quaternion extensor is called
$G$-isometric if, for each $s\in I$, the immersion
$F(s)\otimes G:M^{n-1}\rightarrow {\bf {H^{m}}}:p\mapsto F(s)\otimes G(p)$ is isometric.\\
\begin{lem}\label{Lemma3.1} Let $G:M^{n-1}\rightarrow \bf E^{m}$ be an isometric
immersion of a Riemannian (n-1)-manifold into Euclidean $m$-space
$\bf E^{m}$ and $F:I\rightarrow \bf H$ a unit speed curve in the quaternion plane.
Then\\
(1) the quaternion extensor $\phi=F\otimes G$ is $F$-isometric
if and only if $G$ is unit spherical,\\
(2) the quaternion extensor $\phi=F\otimes G$ is $G$-isometric
if and only if $F$ is unit spherical,\\
(3) the quaternion extensor $\phi=F\otimes G$ is totally real if
and only if either $G$ is spherical or $F(s)=cf(s)$ for some
constant $c\in {\bf H}$ and real-valued function $f$.
\end{lem}
\pf The statements (1) and (2) come from straightforward computations.\\
By a direct computation, the quaternion extensor is totally real
if and only if, for any $s\in I,p\in M^{n-1}$ and $Y\in T_{p}M^{n-1}$, we have
$$Real(\varphi F(s)\bar{F'}(s))\cdot <G(p),Y>=0,$$
where $\bar{F'}$ denotes the quaternionic conjugate of $F'$ and
$Real(\varphi F(s)\bar{F'}(s))$ the real part of $\varphi
F(s)\bar{F'}(s)$ for $\varphi =i,j$ or $k$. Therefore, we have
either $<G(p),Y>=0$ for all $p\in M^{n-1}$, $Y\in T_{p}M^{n-1}$ or
$Real(iF(s)\bar{F'}(s))=
Real(jF(s)\bar{F'}(s))=Real(kF(s)\bar{F'}(s))=0$ for all $s\in I$.
If the first case occurs, then $G$ is spherical. If $F$ is given
by $F(s)=a(s)+ib(s)+jc(s)+kd(s)$, where $a,b,c,$ and $d$ are real
valued functions, and the second condition is true, then we have
the following system of ODEs:
$$\left.\begin{array}{ccc}
ab'-a'b+cd'-c'd & = &0\\
ac'-a'c+b'd-bd' & = &0\\
ad'-a'd+bc'-b'c & = &0
\end{array}\right.$$
By solving this system, we can find that $F(s)=ca(s)$ for a
constant $c\in {\bf{H}}$. \\
\vskip1pc
A submanifold $M^{n-1}$ of $E^{m}$ is said to be of essential
codimension one if locally $M^{n-1}$ is contained in an affine
$n$-subspace of $E^{m}$.\\
\begin{pro}{\label{Proposition 3.2}} Let $G:M^{n-1}\rightarrow \bf E^{m}$
be an isometric immersion of a Riemannian (n-1)-manifold into Euclidean
$m$-space $\bf E^{m}$ and $F:I\rightarrow \bf H$ a unit speed curve.
Then the quaternion extensor $\phi=F\otimes G:I\times M^{n-1}\rightarrow
\bf H^{m}$ is totally geodesic(with respect to the induced metric)
if and only if one of the following two cases occurs:\\
(1) $G:M^{n-1}\rightarrow \bf E^{m}$ is of essential codimension one and
$F(s)=(s+a)c$ for some real number $a$ and some unit quaternion number $c$.\\
(2) $n=2$ and $G$ is a line in $\bf E^{m}$.
\end{pro}
\pf Since $\phi$ is totally geodesic, $\phi_{ss},YZ\phi,Y\phi_{s}$
are tangent vector fields for $Y,Z$ vector fields tangent to the
second component of $I\times M^{n-1}$. By using the fact
$F''(s)\otimes \xi$ is normal to $I\times M^{n-1}$ in $\bf H^{m}$
(via $\phi$) for any unit normal vector field $\xi$ of $M^{n-1}$
in $\bf E^{m}$, we get the following two equations.
\begin{eqnarray}
<F''(s),F''(s)><\xi,G(p)>  = 0  \label{eq:3}\\
<F''(s),F(s)><\xi,h_{G}(Y,Z)>  = 0,\label{eq:4}
\end{eqnarray}
for any vector fields $Y,Z$ tangent to $M^{n-1}$ and for any $s\in I$
and point $p\in M^{n-1}$, where $h_{G}$ is the second fundamental
form of $G:M^{n-1}\rightarrow \bf E^{m}$. \\
We can divide our case as follows:\\
\underline{Case(1)} $F''=0$\\
This case follows from the case(i) in proposition 2.2 in \cite{c97}
so that we can deduce statement (1).\\
\underline{Case(2)} $F''\neq 0$\\
By (\ref{eq:3}), we get $<\xi,G(p)>=0$ for any normal vector field
$\xi$ to $M^{n-1}$ in $\bf E^{m}$ and any point $p\in M^{n-1}$.
Since $\phi$ is totally geodesic, $YZ\phi$ is a tangent vector
field for $Y,Z$ tangent to $M^{n-1}$ in $\bf E^{m}$ which yields
\begin{equation}
0=<F',F><\xi,h_{G}(Y,Z)> \label{eq:5}
\end{equation}
Suppose $G$ is non totally geodesic. Then (\ref{eq:5}) gives $<F',F>=0$
and thus $||F||^{2}$ is a constant. Also, we get
$<F'',F>=0$ because of (\ref{eq:4}). Combining these conditions for
$F$ implies $F'=0$ which is impossible.
Therefore, $G$ must be totoally geodesic.
Since $\phi$ is totally geodesic, $\phi_{ss}$ is tangent to
$I\times M^{n-1}$ in $\bf H^{m}$ so that there exists a tangent
vector $Y$ to $M^{n-1}$ and two real-valued functions $\alpha,\beta$ such that
$$(F''(s)-\alpha(s)F'(s))\otimes G=\beta(s)F(s)\otimes Y.$$
If $F''(s)-\alpha(s)F'(s)=0$ for each $s\in I$, then we have
$F''=0$ which is the contradiction to our case. Thus,
$G(p)$ is tangent to a vector $Y$ which implies that $n=2$.
Furthermore, $G$ is a line in $\bf E^{m}$ since $G$ is totally geodesic.
The converse can be proved easily.
\begin{pro}{\label{proposition 3.3}} let $\imath:S^{n-1}\rightarrow \bf E^{n}$
be the inclusion of the unit hypersphere of $\bf E^{n}$(centered at the origin).
Then every quaternion extensor of $\imath$ via a unit speed curve $F$ in $\bf H$
is a Lagrangian $H$-umbilical submanifold of $\bf H^{n}$ unless $F(s)=(s+a)c$
for some real number $a$ and some unit quaternion number $c$.
\end{pro}
\pf By a direct computation, we can easily see that $\phi=F\otimes \imath$ is a
Lagrangian $H$-umbilical submanifold satisfying
\begin{eqnarray*}
h(e_{1},e_{1}) & = & \lambda_{I} Ie_{1}+\lambda_{J} Je_{1}+\lambda_{K} Ke_{1},   \\
h(e_{1},e_{j}) & = & \mu_{I} Ie_{j}+\mu_{J} Je_{j}+\mu_{K} Ke_{j}, \mbox{   for } j=2,...,n   \\
h(e_{j},e_{j}) & = & \mu_{I} Ie_{1}+\mu_{J} Je_{1}+\mu_{K} Ke_{1}, \mbox{   for } j=2,...,n   \\
h(e_{j},e_{k}) & = & 0, \mbox{           for } j\neq k =2,...,n,
\end{eqnarray*}
where $\lambda_{\varphi}=<F'',\varphi F'>$ and
$\mu_{\varphi}=<(\frac{F}{\parallel F \parallel})^{'},
\varphi(\frac{F}{\parallel F \parallel})>$ for $\varphi=I,J$ or
$K$ and $\{e_{1},e_{2},...,e_{n}\}$ is an orthonormal local frame
field. Without difficulty, we can get that $\phi$ is totally
geodesic if $F(s)=(s+a)c$ for some real number $a$ and some unit
quaternion number $c$.
\section{Main theorem}
The main result of this section is to classify Lagrangian $H$-umbilical submanifolds
of quaternion Euclidean space. To do this, we need to review the Lagrangian
pseudo-sphere in $\bf C^{n}$(\cite {c97}). For a real number $b>0$, let
$F:\bf R \rightarrow \bf C$ be the unit speed curve given by
\[
F(s)=\displaystyle{\frac{e^{2bsi}+1}{2bi}}
\]
With respect to the induced metric, the complex extensor $\phi =F \otimes \imath$ of
the unit hypersphere of $E^{n}$ via $F$ is a Lagrangian isometric immersion of an
open portion of an n-sphere $S^{n}(b^{2})$ of sectional curvature $b^{2}$ into
$\bf C^{n}$. It is called a Lagrangian pseudo-sphere. It has been shown that
it is a Lagrangian $H$-umbilical submanifold in $\bf C^{n}$ satisfying the
following second fundamental form:
\begin{eqnarray}\label{eq:10}
h(e_{1},e_{1}) & = & 2bJe_{1}, \hskip1pc h(e_{i},e_{i})  =  bJe_{1},\hskip.3pc i\geq 2 \\
h(e_{1},e_{j})&  = & bJe_{j}, \hskip1pc h(e_{j},e_{k})  =
0,\hskip.3pc \mbox{for } j\neq k=2,...,n,\nonumber
\end{eqnarray}
for some nontrivial function $b$ with respect to some suitable orthonormal
local frame field.
Up to rigid motions in $\bf C^{n}$, it is unique.
\begin{thm}\label{theorem:4.1} Let $n \geq 3$ and $L:M\rightarrow \bf H^{n}$
be a Lagrangian $H$-umbilical isometric immersion.\\
We have one of these three cases:\\
(A) $M$ is flat or, \\
(B) up to rigid motions of $\bf H^{n}$, $L$ is a Lagrangian pseudo-sphere
in $\bf C^{n}$, or\\
(B) up to rigid motions of $\bf H^{n}$, $L$ is a quaternion extensor of
the unit hypersphere of $\bf E^{n}$.
\end{thm}
\pf Let $n\geq 3$ and $L:M\rightarrow \bf H^{n}$ be a Lagrangian
$H$-umbilical isometric immersion whose second fundamental form is given by
\begin{eqnarray}
h(e_{1},e_{1})& = & \lambda_{1} Ie_{1}+\lambda_{2} Je_{1}+\lambda_{3} Ke_{1},\label{eq:6}\\
h(e_{1},e_{j})& = & \mu_{1} Ie_{j}+\mu_{2} Je_{j}+\mu_{3} Ke_{j}, \mbox{   for } j=2,...,n \nonumber  \\
h(e_{j},e_{j})& = & \mu_{1} Ie_{1}+\mu_{2} Je_{1}+\mu_{3} Ke_{1}, \mbox{   for } j=2,...,n  \nonumber \\
h(e_{j},e_{k})& = & 0, \mbox{           for } j\neq k =2,...,n,\nonumber
\end{eqnarray}
for some functions $\lambda_{i},\mu_{i}$(i=1,2,3) with respect to some suitable
orthonormal local frame fields $\{e_{1},e_{2},...,e_{n}\}$ with the dual 1-forms
$\omega_{1},...,\omega_{n}$.
Let $(\omega_{A}^{B})$, $A,B=1,...,n$ be the connection form on $M$ definded by
$\omega_{i}^{j}(e_{k})=<\widetilde{\nabla}_{e_{k}}e_{i},e_{j}>$ for $i,j,k=1,...,n$.
By (\ref{eq:6}) and Codazzi equation, we have
\begin{eqnarray}\label{eq:7}
e_{1}(\mu_{1})=(\lambda_{1}-2\mu_{1})\omega_{1}^{j}(e_{j})+\lambda_{2}\mu_{3}-\lambda_{3}\mu_{2} \\
e_{1}(\mu_{2})=(\lambda_{2}-2\mu_{2})\omega_{1}^{j}(e_{j})+\lambda_{3}\mu_{1}-\lambda_{1}\mu_{3}\nonumber \\
e_{1}(\mu_{3})=(\lambda_{3}-2\mu_{3})\omega_{1}^{j}(e_{j})+\lambda_{1}\mu_{2}-\lambda_{2}\mu_{1}\nonumber
\end{eqnarray}
\begin{eqnarray}\label{eq:8}
e_{j}(\lambda_{1})&=&(\lambda_{1}-2\mu_{1})\omega_{1}^{j}(e_{1}) \\
e_{j}(\lambda_{2})&=&(\lambda_{2}-2\mu_{2})\omega_{1}^{j}(e_{1})\nonumber \\
e_{j}(\lambda_{3})&=&(\lambda_{3}-2\mu_{3})\omega_{1}^{j}(e_{1}) \mbox{ for } j=2,...,n \nonumber
\end{eqnarray}
\begin{eqnarray}\label{eq:9}
(\lambda_{1}-2\mu_{1})\omega_{1}^{k}(e_{j})&=&0 \\
(\lambda_{2}-2\mu_{2})\omega_{1}^{k}(e_{j})&=&0 \nonumber \\
(\lambda_{3}-2\mu_{3})\omega_{1}^{k}(e_{j})&=&0 \mbox{ for } k\neq j=2,...,n \nonumber
\end{eqnarray}
\begin{eqnarray}\label{eq:10}
e_{j}(\mu_{1})&=&3\mu_{1}\omega_{1}^{j}(e_{1})\\
e_{j}(\mu_{2})&=&3\mu_{2}\omega_{1}^{j}(e_{1})\nonumber \\
e_{j}(\mu_{3})&=&3\mu_{3}\omega_{1}^{j}(e_{1}) \mbox{ for } j=2,...,n \nonumber
\end{eqnarray}
\begin{equation}\label{eq:11}
\mu_{1}\omega_{1}^{k}(e_{1})=\mu_{2}\omega_{1}^{k}(e_{1})=\mu_{3}\omega_{1}^{k}(e_{1})=0, k\neq j=2,...,n
\end{equation}
Here, (\ref{eq:9}) and (\ref{eq:11}) hold only for $n\geq 3$.
\vskip.1pc
We can divide our case as follows.\\
(A) If $M$ is of constant sectional curvature, then (\ref{eq:6}) implies that
$\mu_{1}(\lambda_{1}-2\mu_{1})+\mu_{2}(\lambda_{2}-2\mu_{2})+\mu_{3}(\lambda_{3}-2\mu_{3})=0$.
If $\mu_{1}=\mu_{2}=\mu_{3}=0$, then $M$ is flat.
\vskip.2pc
From now on, we assume that there exists at least one $\mu_{i}$ which is not
identically zero.
Furthermore, the equations in (\ref{eq:9}) provide the following two cases:\\
(a) $\lambda_{i}=2\mu_{i}$ for $i=1,2,3$\\
This condition and the equations in (\ref{eq:7}) and (\ref{eq:8})
imply that $\mu_{1},\mu_{2}$ and $\mu_{3}$ are constants.
Moreover, by Gauss equation, $M$ is a real-space form of constant
sectional curvature $\mu_{1}^{2}+ \mu_{2}^{2}+\mu_{3}^{2}$, say
$b^{2}\neq 0$. By making a proper translation and rescaling, we
can say that $M$ satisfies the second fundamental form given in
(\ref{eq:10}). Moreover, we also can check that the first normal
space is parallel with respect to the normal connection so that by
applying the result of Erbacher \cite{erbacher}, $M$ can be
immersed into complex Euclidean space $\bf C^{n}$ which implies
that $M$ is a Lagrangian pseudo-sphere in $\bf C^{n}$.
\\[0.1in]
(b) There exists one $i$ such that $\lambda_{i}\neq 2\mu_{i}$, and also
we have $\mu_{1}(\lambda_{1}-2\mu_{1})+\mu_{2}(\lambda_{2}-2\mu_{2})+\mu_{3}(\lambda_{3}-2\mu_{3})=0$.\\
By (\ref{eq:9}), we know that $\omega _{1}^{j}(e_{k})=0$ for $k\neq j=2,...,n$.\\
By (\ref{eq:7}) and (\ref{eq:11}), we get
\begin{equation}\label{eq:12}
\omega_{1}^{j}=\displaystyle{\frac{e_{1}(\mu_{2})-\lambda_{3}\mu_{1}+\lambda_{1}\mu_{3}}
{\lambda_{2}-2\mu_{2}}\omega^{j}}
\end{equation}
Suppose we define $f=\displaystyle{\frac{e_{1}(\mu_{2})-\lambda_{3}\mu_{1}+
\lambda_{1}\mu_{3}}{\lambda_{2}-2\mu_{2}}}$.
Let $D$ be the distribution spanned by $e_{1}$ and $D^{\perp}$ be the
distribution spanned by $\{e_{2},e_{3},...,e_{n}\}$. Since
$\omega _{1}^{j}(e_{k})=0$ for $k\neq j=2,...,n$, the
distribution $D^{\perp}$ is integrable. Also, the distribution
$D$ is intergrable since it is 1-dimensional. Therefore there
exists local coordinates $\{x_{1},x_{2},...,x_{n}\}$ such that
$e_{1}=\frac{\partial}{\partial x_{1}}$ and $D^{\perp}$ is
spanned by $\{\frac{\partial}{\partial x_{2}},...,\frac{\partial}
{\partial x_{n}}\}$. Using (\ref{eq:12}), we obtain
\begin{equation}\label{eq:13}
<\nabla_X Y,e_{1}>=-f<X,Y>
\end{equation}
It implies that $D^{\perp}$ is a spherical distribution and
futhermore, each leaf of $D^{\perp}$ is of constant sectional
curvature $\mu_{1}^{2}+\mu_{2}^{2}+ \mu_{3}^{2}+k^{2}$. Now, by
applying a result of Hiepko \cite{h}, $M$ is isometric to a warped
product $I \times_{\omega(s)} S^{n-1}$, where
$S^{n-1}$ is the unit $(n-1)$ sphere and $\omega(s)$ is a warping function.\\
Using the spherical coordinates $\{u_{2},...,u_{n}\}$ on the unit sphere,
we can choose the metric
\[g=ds^{2}+\omega^{2}(s)\{du_{2}^{2}+ \cos^{2}u_{2}du^{2}_{3}+
\cdots + \cos^{2} u_{2} \cdots \cos^{2} u_{n-1}du_{n}^{2}\}\]
on $I \times_{\omega(s)} S^{n-1}$. By using this metric $g$, we have
\begin{equation}\begin{array}{l}\label{eq:14}
\nabla_{\frac{\partial}{\partial s}}\frac{\partial}{\partial s}=0,
\mbox{    }\nabla_{\frac{\partial}{\partial s}}{\frac{\partial}{\partial u_{k}}=
\frac{\omega'}{\omega}\frac{\partial}{\partial u_{k}},
\mbox{   }\nabla_{\frac{\partial}{\partial u_{2}}}\frac{\partial}{\partial u_{2}}=
-\omega\omega' \frac{\partial}{\partial s}},\\
\nabla_{\frac{\partial}{\partial u_{i}}}\frac{\partial}{\partial u_{j}}=-\tan u_{i}
\frac{\partial}{\partial u_{j}}, 2\leq i<j,\\
\nabla_{\frac{\partial}{\partial u{j}}}\frac{\partial}{\partial u_{j}}=
-\omega\omega'\cos^{2}u_{2}\cdots \cos^{2}u_{j-1}\frac{\partial}{\partial s}\\
          +\stackrel{n-1}{\sum}_{k=2}\frac{\sin 2u_{k}}{2}\cos^{2}u_{k+1}
          \cdots \cos^{2}u_{j-1}\frac{\partial}{\partial u_{k}}, j\geq 2
\end{array}
\end{equation}
By substituting $X=Y=\frac{\partial}{\partial u_{2}}$ into (\ref{eq:13}) and
using (\ref{eq:14}), we get
$$\frac{\omega'}{\omega}=f$$
Furthermore,  computing the sectional curvatures spanned by
$\frac{\partial}{\partial u_{2}},\frac{\partial}{\partial u_{3}}$
and $\frac{\partial}{\partial s},\frac{\partial}{\partial u_{2}}$
derive the following condition for the warping fuction:
$$\frac{1}{\omega^{2}}-f^{2}=-\frac{\omega''}{\omega}=
\mu_{1}^{2}+\mu_{2}^{2}+\mu_{3}^{2}=\bar{\mu}^{2}$$
Note here that $\bar{\mu}$ is a constant by our assumption.
These conditions provide a differential equation
\[f^{2}+f'+\bar{\mu}^{2}=0\]
and then
\[\omega(s)=\cos \bar{\mu}s,f(s)=-\bar{\mu}\tan \bar{\mu}s\]
By doing the same procedure in theorem 4.1 \cite{c97}, we
obtain that $L$ is a quaternion extensor of the unit hypersphere of $\bf E^{n}$.
\vskip.3pc
(B) Now, we assume that $M$ does not contain open subset of constant
sectional curvature. Then
$$U=:\{p\in M:\mu_{1}(\lambda_{1}-2\mu_{1})+\mu_{2}(\lambda_{2}-
2\mu_{2})+\mu_{3}(\lambda_{3}-2\mu_{3})\neq 0 \mbox{ at }p\}$$
is an open dense subset of $M$.\\
By (\ref{eq:9}), on $U$, we obtain
$$\omega_{1}^{j}(e_{j})=\displaystyle{\frac{\mu_{1}e_{1}(\mu_{1})+
\mu_{2}e_{1}(\mu_{2})+\mu_{3}e_{1}(\mu_{3})}{\mu_{1}(\lambda_{1}-
2\mu_{1})+\mu_{2}(\lambda_{2}-2\mu_{2})+\mu_{3}(\lambda_{3}-2\mu_{3})}}$$
Since $n\geq 3$, $\omega_{1}^{j}(e_{k})=0$ for $k\neq j=2,...,n$ on $U$.
Therefore,
$$\omega_{1}^{j}=\displaystyle{\frac{\mu_{1}e_{1}(\mu_{1})+
\mu_{2}e_{1}(\mu_{2})+\mu_{3}e_{1}(\mu_{3})}{\mu_{1}(\lambda_{1}-2\mu_{1})+
\mu_{2}(\lambda_{2}-2\mu_{2})+\mu_{3}(\lambda_{3}-2\mu_{3})}}\omega^{j}$$
Let's define $\bar {f}=\displaystyle{\frac{\mu_{1}e_{1}(\mu_{1})+
\mu_{2}e_{1}(\mu_{2})+\mu_{3}e_{1}(\mu_{3})}{\mu_{1}(\lambda_{1}-
2\mu_{1})+\mu_{2}(\lambda_{2}-2\mu_{2})+\mu_{3}(\lambda_{3}-2\mu_{3})}}$
Using this $\bar{f}$, the conclusion follows in the same way we have seen the above case.
\vskip1pc
Similar to the complex case, we have the following same result
for Lagrangian $H$-umbilical surface in quaternion Euclidean space $\bf H^{2}$.
\begin{thm}\label{theorem:4.2} Let $L:M\rightarrow \bf H^{2}$ be a
Lagrangian $H$-umbilical surface satisfying
\begin{eqnarray*}
h(e_{1},e_{1})=\lambda_{1}Ie_{1}+\lambda_{2}Je_{1}+\lambda_{3}Ke_{1},\\
h(e_{2},e_{2})=\mu_{1}Ie_{1}+\mu_{2}Je_{1}+\mu_{3}Ke_{1},\\
h(e_{1},e_{2})=\mu_{1}Ie_{2}+\mu_{2}Je_{2}+\mu_{3}Ke_{2},
\end{eqnarray*}
such that the integral curves of $e_{1}$ are geodesics in $M$. Then\\
(1) $M$ is flat or, \\
(2) up to rigid motions of $\bf H^{2}$, $L$ is a Lagrangian
pseudo-sphere in $\bf C^{2}$ or \\
(3) up to rigid motions of $\bf H^{2}$, $L$ is a quaternion
extensor of the unit circle of $\bf E^{2}$.
\end{thm}
\pf The proof is similar to theorem 4.1 with some modifications.
\vskip.5pc
\underline{Remark} The explicit description of flat Lagrangian
$H$-umbilical submanifolds in a quaternion Euclidean space will be
discussed in \cite{oh3}.
\\[0.1in]
The author would like to thank Prof. B.Y. Chen for suggesting the
problem and for useful discussions on this topic.
\vskip.5pc
Yun Myung Oh\\
Department of Mathematics\\
Michigan State University\\
E. Lansing, MI 48824\\
U.S.A.\\
\end{document}